# Blow-Up Results of Time Fractional Heat Equation With Nonlinear Neumann Boundary Conditions

**Hind Ghazi Hameed[1,2], Burhan Selcuk [3] and Maan A. Rasheed[4]***

[1] Postgraduate education institute, Department of Mathematics, Karabuk University, Turkey -
[2] Department of Mathematics, College of Education for Pure Science, Tikrit University, Iraq
[3] Department of Computer Science, Faculty of Science, Ankara University, Turkey
[4] Department of Mathematics, College of Science for Women, University of Baghdad, Iraq

[1] 2038179003@karabuk.edu.tr, hind.g.hameed@tu.edu.iq
[2] burhanselcuk@ankara.edu.tr

*Corresponding Author: maan.a@csw.uobaghdad.edu.iq

**Abstract:** The study of blow-up solution of time-fractional heat equations is of significant and wide-ranging interest for its multitude of applications. These types of equations are used to model several real problems in science and engineering. This article is concerned with the blow-up solutions of a time fractional heat equation subject to nonlinear Neumann boundary conditions of power type. Firstly, the local global exitance of positive solutions and blows up in finite time are studied, under some restricted conditions . Secondly, the blow-up set is investgated showing that the blow up can only occur at a boundary point.

**1. Introduction**

we consider the initial- boundary value problem of time-Caputo fractional heat equation:

$$\begin{cases} {}^c D_t^\alpha u(x,t) = u_{xx}(x,t), & 0 < x < 1, 0 < t < T, \\ u_x(0,t) = 0, u_x(1,t) = u^p(1,t), & 0 < t < T, \\ u(x,0) = u_0(x), & 0 \leq x \leq 1. \end{cases} \qquad (1)$$

where ${}^c D_t^\alpha u$ denotes Caputo fractional derivative operator, defined as

$$ {}^c D_t^\alpha u(x,t) = \frac{1}{\Gamma(1-\alpha)} \int_0^t (t-\tau)^{-\alpha} \frac{\partial u(x,\tau)}{\partial \tau} d\tau,$$

and $\Gamma$ is Gamma function with $0 < \alpha < 1$. Also, we assume that $p > 1$, and $T$ is the maximal existence time. The initial function $u_0(x)$ satisfies the compatibility conditions:

$$u_0'(0) = 0, u_0'(1) = u_0^p(1)$$

Blow-up is a phenomenon that refers to a case where a solution to a partial differential equation becomes unbounded or infinite within a finite time interval. In some Latin-derived languages, as an explosion. We can generalize the concept of blowing up more broadly, as it represents the phenomenon in which solutions cannot be globally continued over time, resulting from the infinite growth of the nonlinear terms that show the evolutionary process.

The heat equation or diffusion equation describes the physical processes that involve the transfer of a given quantity (such as heat or matter) from a high concentration area to a low concentration area due to diffusion, a physical phenomenon that appears to be widespread. Diffusion results from Brownian motion, which is defined as the random distribution of atoms, particles, and molecules floating or suspended in a fluid, resulting in complete mixing. At the macro-level, it is defined by the famous Fick's laws that describe diffusion or Fourier's laws, also known as the law of heat conduction [1]. This classical heat equation is linear and given by $u_t = c^2 \cdot u_{xx}$ where c is diffusion coefficient and u is temperature or concentration. It can be solved easily. However, the problem becomes more complicated when a source is added to the equation, as it becomes a semi-linear heat equation.

In physical problems [2], the time-fractional derivative presents the nonlocal nature. Therefore, the time Caputo time-fractional heat equation subject to Neumann boundary conditions can be considered as a partial differential equation (PDE) that incorporates fractional calculus in modeling the diffusion of heat over time. It is an extension of the classical heat equation, where the presence of the Caputo fractional derivative in the time variable introduces memory effects and nonlocal behavior. These types of problems can be used for modelling several real-life problems involving memory, delay effects and non-local descriptions.

While classical diffusion equations have been discussed and has had results in literature. However, for fractional diffusion equations, the answer to this question is still in its infancy. For several reasons, the existing literature refers to time-fractional diffusion equations with relatively little effort devoted to studying the abstract mathematical properties of fractional differential equations because attention has been paid to physical modeling.

The study of blow-up phenomena first appeared in the 1940s and 1950s. More comprehensive studies on this subject were conducted in the 1960s, especially by [3-6]. Many researchers have studied the behavior and blow-up of solutions of linear and semilinear heat equations with nonlinear boundary conditions (Neumann). In particular, the blow-up properties of the following two problems have been studied in detail by many researchers [7-13]:

$$\begin{cases} u_t(x,t) = \Delta u(x,t) + u^p, & \text{in } \Omega \times (0,T), \\ \dfrac{\partial u}{\partial \eta} = 0, & \text{on } \partial\Omega \times (0,T), \\ u(x,0) = u_0(x), & \text{in } \Omega. \end{cases} \qquad (2)$$

$$\begin{cases} u_t(x,t) = \Delta u(x,t), & \text{in } \Omega \times (0,T), \\ \dfrac{\partial u}{\partial \eta} = u^p, & \text{on } \partial\Omega \times (0,T), \\ u(x,0) = u_0(x), & \text{in } \Omega. \end{cases} \qquad (3)$$

Lin and Wang [14] focus on studying the blow up properties of the one-dimensional semi-linear heat equation with Neumann boundary conditions;

$$\begin{cases} u_t(x,t) = u_{xx}(x,t) + u^p(x,t), & 0 < x < 1, t > 0 \\ u_x(0,t) = 0, u_x(1,t) = u^q(1,t), & t > 0 \\ u(x,0) = u_0(x), & 0 \leq x \leq 1 \end{cases} \qquad (4)$$

The study clarified the interaction between the two nonlinear terms and the role of each in determining the properties, rate, and location of the blow-up. They used analytical methods, including the comparison principle, energy methods, and integral equations, to reach the results. They proved that the blow-up occurs in a finite time when $max(p,q) > 1$. They reach the blow-up rates and concluded that the blow-up occurs at boundary only when $x = 1$.

The concepts of blow-up and quenching can be transformed into each other by means of a transformation [15]. Ozalp and Selcuk studied the behavior of solutions and determined standards for blow-up and quenching to these equations. Later, they considerd heat equation with nonlinear boundary conditions:

$$\begin{cases} u_t(x,t) = u_{xx}(x,t), & 0 < x < 1, t > 0 \\ u_x(0,t) = u^p(0,t), u_x(1,t) = u^q(1,t), & t > 0 \\ u(x,0) = u_0(x), & 0 \leq x \leq 1 \end{cases} \quad (5)$$

They studied this problem using both analytical and numerical tools. They developed standards for blowing up and putting out an blow up based on how similar the two phenomena used to be. The study results include criteria for blowing up and cooling in the heat equation and the non-linear parabolic equation with non-linear boundary conditions. The authors also showed that blowing up and quenching are equivalent in these equations, and they used numerical models to back up their claims.

In the article [16], Levine focused on studying the phenomenon of quenching in solutions of linear parabolic equations with nonlinear boundary conditions, where the study focused on determining the conditions under which the solutions reach a critical value, where this condition is called quenching.

$$\begin{cases} u_t(x,t) = u_{xx}(x,t), & 0 < x < L, \quad t > 0 \\ u(0,t) = 0, u_x(L,t) = \phi(u(L,t)), & t > 0 \\ u(x,0) = u_0(x), & 0 \leq x \leq L \end{cases} \quad (6)$$

He used analytical methods, including the principles of maximum values, energy methods, and Green's function, to reach the results. The uniqueness of the solutions was proven and the conditions that allow their extension after quenching were extended. He also proved that quenching occurs in a finite time when $L_0 > L$ at the boundary point $x = L$ and $u_t(L,t)$ becomes unbounded, and concluded that quenching does not occur and the solutions remain globally when $L_0 < L$.

Many blow-up and quenching problems with different initial and boundary conditions in different types of equations have been considered [17-22] In addition, many blow up problems have recently been added to the literature by using fractional derivatives instead of classical derivatives [23-32].

Recently, studies have been focused on the behavior of blow-up solutions of time-fractional diffusion equations. For instance, Subedi and Vatsala have studied in [33] some blow up properties for time fractional one-dimensional semilinear reaction-diffusion equation subject to homogenuous Dirchlet boundary conditions.

$$\begin{cases} {}^cD_t^\alpha u(x,t) = u_{xx}(x,t) + u^p(x,t), & 0 < x < 1, \quad 0 < t < T \\ u(0,t) = 0, u(L,t) = 0, & 0 < t < T \\ u(x,0) = 0, & 0 \leq x \leq 1 \end{cases} \quad (7)$$

They used analytical methods, including comparison principles methods where solutions of fractional equations are compared with solutions of integer equations where the blow up behavior of equations is well known. They also used the construction of the lower solution, where lower solutions of fractional equations are constructed using solutions of integer equations. They also used Green's function to achieve the results. They concluded that for ordinary fractional equations, the lower solutions extracted from integer equations blow up in a finite time. These lower solutions provide upper limits on the time in which the blow up of fractional equations occurs. As for fractional diffusion reaction equations, if solutions of integer equations blow up, the solutions of fractional equations will also explode under suitable conditions. Blow up time in the fractional case depends on the blow-up time in the integer equations and is modified using the fractional derivative. To summarize, the researchers showed that the blow-up behavior of diffusion reaction equations is very similar to that of integer equations when similar conditions are met.

Motivated by studies [14] and [33] in this paper, it is works on blow up of timefractional heat equations with a nonlinear Neumann boundary condition. Here, it is extended the classical heat equation to the fractional heat equation and took a nonlinear Neumann boundary condition instead of the Dirichlet boundary condition.

This paper is divided into five sections. In section two, we recall some basic preliminaries. In the third section, the conditions of global existence and blow-up in finite time are obtained. Section four is devoted to studying the blow-up point. Finally, in the last section, we give the main conclusions, theoretical and practical contributions, and possible directions for future studies.

## 2. Preliminaries

In this section, we recall some basic definitions and Auxiliary lemma that we shall need to prove the main results.

**Definition 2.1.** [33] A solution $u(x,t) \in C^{2,\alpha}([0,1] \times [0,T))$ of the equation (1) blows up in a finite-time, i.e, the existence of a $T = T(u_0) < \infty$ such that

$$\lim_{t \to T^-} max\{u(x,t): 0 \leq x \leq 1\} \to \infty$$

The blow up time of the equation (1) is denoted as $T$.

**Definition 2.2.** [33] A function $v(x,t) \in C^{2,\alpha}([0,1] \times [0,T))$ is called to be lower solution of the equation (1), if

$$\begin{cases} {}^C D_t^\alpha v(x,t) - v_{xx}(x,t) \leq 0, & 0 < x < 1, 0 < t < T \\ v_x(0,t) \leq 0, v_x(a,t) \leq v^p(1,t), & 0 < t < T \\ v(x,0) \leq 0, & 0 \leq x \leq 1 \end{cases} \quad (8)$$

and $w(x,t) \in C^{2,\alpha}([0,1] \times [0,T))$ is called the upper solution of the equation (1), if the reverses of the inequalities in (8) are held.

**Definition 2.3.** [34] For two parameters, $\alpha, r > 0$, the Mittag Leffler function is defined as follows:

$$E_{\alpha,r}(z) = \sum_{k=0}^{\infty} \frac{z^k}{\Gamma(\alpha k + r)}$$

Specially,

$$E_{\alpha,r}(\lambda(t-t_0)^\alpha) = \sum_{k=0}^{\infty} \frac{(\lambda(t-t_0)^\alpha)^k}{\Gamma(\alpha k + r)}$$

Also, for $t_0 = 0$ and $r = 1$, we get,

$$E_{\alpha,1}(\lambda t^\alpha) = \sum_{k=0}^{\infty} \frac{(\lambda t^\alpha)^k}{\Gamma(\alpha k + 1)}$$

and for $t_0 = 0$ and $r = \alpha$, we get,

$$E_{\alpha,\alpha}(\lambda t^\alpha) = \sum_{k=0}^{\infty} \frac{(\lambda t^\alpha)^k}{\Gamma(\alpha(k+1))}$$

where $0 < \alpha < 1$.

Further if $\alpha = 1$, then $E_{\alpha,\alpha}(\lambda t^\alpha) = E_{\alpha,1}(\lambda t^\alpha) = e^{\lambda t}$.

**Remark 2.1**

It is known that the Maximum principles and Hopf's lemma are considered essential tools to prove several qualitative properties of classical diffusion equations, including blow-up and quenching problems. In this article, the Maximum principles and Hopf's lemma for time-fractional equations, obtained by some researchers, in [35-37], are employed to get the main results.

**Auxiliary Lemma**

By following the logical process of Levine [16] applied for the quenching problem, the following auxiliary theorem for the solution of problem (1), $u(x,t)$, can be easily obtained using the fractional maximum principles and the fractional Hopf's lemma.

**Lemma 3.1.** If $u$ solves (1) in $(0,1) \times (0,T)$, then:

- $u > 0$ in $(0,1) \times (0,T)$.
- $u_x(x,t) \geq 0$ in $(0,1) \times (0,T)$ where $u_x(x,0) \geq 0$ in $(0,1)$.
- ${}^cD_t^\alpha(u(x,t)) > 0$ in $(0,1) \times (0,T)$, where $u_{xx}(x,0) \geq 0$ in $(0,1)$.

**3. Blow-up and Global Exsitance**

This section is devoted to studying the conditions of global existence and blow-up of problem (1).

**Theorem 3.1**. A solution of problem (1) exists globally, if $p \leq 1$.

Proof. The sufficient condition must be proven first. We need to show that the solution exists globally. The auxiliary function is $v(x,t) = CE(Kt^\alpha)e^{Lx^2}$ where C, K and L sufficiently large constants to be determined later and $E(Kt^\alpha)$ is MittagLeffler function ${}^cD_t^\alpha(E(Kt^\alpha)) = KE(Kt^\alpha)$. It is calculated derivation of $v$ in the following

$$^cD_t^\alpha(v(x,t)) = CKE(Kt^\alpha)e^{Lx^2}$$
$$v_x(x,t) = 2CLxE(Kt^\alpha)e^{Lx^2}$$
$$v_{xx}(x,t) = (4L^2x^2 + 2L)CE(Kt^\alpha)e^{Lx^2}$$

and therefore, we get;

$$^cD_t^\alpha(v(x,t)) - v_{xx}(x,t) \geq 0$$

where $K = 4L^2 + 2L$.

Apply boundary conditions:

$$v_x(0,t) = 0$$
$$v_x(1,t) = 2CLE(Kt^\alpha)e^L \geq (2CLE(Kt^\alpha)e^L)^p = v^p(1,t)$$

where $p \leq 1$. Also, $v(x,0) \geq 0$ on the initial line. It is seen that $v(x,t) = CE(Kt^\alpha)e^{Lx^2}$ is an upper solution of the problem (1). It is obtained from Definition 2.1 and Definition 2.2 that problem (1) exits globally where $p \leq 1$.

**Theorem 3.2.** A solution of the problem (1) blows up in a finite time and an upper limit for the blow-up time is

$$\left(\frac{\Gamma(\alpha+1)}{p-1}\right)^{\frac{1}{\alpha}} \left(\frac{1}{l_0}\right)^{\frac{\alpha}{p-1}}$$

where $p > 1, u_x(x,0) \geq 0$ and and $l_0(>0)$ is the initial function of the following problem

$$l_t = l^p(t)$$

*Proof.* Now, suppose that $p > 1$ and show that the solution of problem (1) blows up in finite time. Let $m(t) = \int_0^1 u(x,t)dx$.

$$^cD_t^\alpha(m(t)) = \int_0^1 {^cD_t^\alpha}u(x,t)dx = \int_0^1 u_{xx}(x,t)dx = u^p(1,t)$$

From Lemma 3.2, $u(1,t) \geq u(x,t)$ where $x \in [0,1)$. Then, from the above equation and given $> 1$ :

$$^cD_t^\alpha(m(t)) \geq m^p(t) \tag{9}$$

Vatsala ve Subedi investigate the solution of the following problem in [33];

$$^cD_t^\alpha(n(t)) = n^p(t), n(0) = m(0) > 0 \tag{10}$$

It is clear that the solution of (9) is an upper solution of (10). Unfortunately, (10) does not have an explicit solution. The problem (10) with classical differentiation is expressed as follows;

$$l_t = l^p(t), l(0) > 0, \tag{11}$$

They show that $l\left(\frac{t^\alpha}{\Gamma(\alpha+1)}\right)$ is a lower solution of the problem (10) where $p > 1$. Also, they obtained that the blow-up time of the problem (11) is

$$t_l = \left(\frac{\Gamma(\alpha+1)}{p-1}\right)^{\frac{1}{\alpha}} \left(\frac{1}{l_0}\right)^{\frac{\alpha}{p-1}}.$$

As a result, a solution of the problem (1) blows up in a finite time. Since the solution to problem (11) is a subsolution of problem (9) and it blows up, the solution to problem (9) also blows up. That is, $u(x,t)$ also blows up. In addition, it is obtained an upper limit of the blow-up time of the problem (9) is

$$\left(\frac{\Gamma(\alpha+1)}{p-1}\right)^{\frac{1}{\alpha}} \left(\frac{1}{l_0}\right)^{\frac{\alpha}{p-1}}.$$

since $t_m \leq t_n \leq t_l$ from Definition 2.2.

For problem (1), the following results can be concluded with the help of Theorem 3.1 and Theorem 3.2;

**Corollary 3.1.**

If $p > 1$, the solution can blow up in a finite time,

If $p \leq 1$, the solution exists globally.

### 4. Blow-up Point

In this section, we consider the blow-up set of problem (1) and show that it can only be achieved at a boundary point.

**Theorem 4.3.** If $u_x(x,0) \geq 0$ and $u_{xx}(x,0) \geq 0$ in $[0,1]$, then $x = 1$ is the only blow-up point for problem (1).

*Proof.* Let $d_1 \in [0,1), d_2 \in [d_1, 1), \tau \in (0,T)$ and $\epsilon > 0$. Define

$$F(x,t) = u_x(x,t) - \epsilon(x-d_1)u^p(x,t) \text{ in } [d_1, d_2] \times [\tau, T]$$

where $p > 1$ and $\epsilon$ is a sufficient small constant. The 1st and 2nd derivatives of the auxiliary function $F(x,t)$ with respect to x are obtained as follows;

$$F_x(x,t) = u_{xx}(x,t) - \epsilon u^p(x,t) - \epsilon p(x-d_1)u^{p-1}(x,t)u_x(x,t),$$
$$F_{xx}(x,t) = u_{xxx}(x,t) - 2\epsilon p u^{p-1}(x,t)u_x(x,t) - \epsilon p(x-d_1)u^{p-1}(x,t)u_{xx}(x,t)$$
$$-\epsilon p(p-1)(x-d_1)u^{p-2}(x,t)u_x^2(x,t).$$

On the other hand:

$$^cD_t^\alpha u^p(x,t) = \frac{p}{\Gamma(1-\alpha)} \int_0^t \frac{u^{p-1}(x,s)u_s(x,s)}{(t-\tau)^\alpha} ds$$

By generalized Leibinz rule, we have

$$^cD_t^\alpha u^p(x,t) = pu^{p-1}(x,t)\,^cD_t^\alpha u(x,t) - \frac{p(p-1)}{\Gamma(1-\alpha)}\int_0^t \frac{u^{p-2}(x,s)(u(x,t)-u(x,s))u_s(x,s)}{(t-\tau)^\alpha}ds.$$

$$^cD_t^\alpha F(x,t) = {}^cD_t^\alpha u_x(x,t) - \epsilon(x-d_1)\,^cD_t^\alpha u^p(x,t)$$
$$= u_{xxx}(x,t) - \epsilon(x-d_1)pu^{p-1}(x,t)\,^cD_t^\alpha u(x,t)$$
$$+\epsilon(x-d_1)\frac{p(p-1)}{\Gamma(1-\alpha)}\int_0^t \frac{u^{p-2}(x,s)(u(x,t)-u(x,s))u_s(x,s)}{(t-\tau)^\alpha}ds$$

Since is increasing t, it follows that

$$^cD_t^\alpha F(x,t) \geq u_{xxx}(x,t) - \epsilon(x-d_1)pu^{p-1}(x,t)\,^cD_t^\alpha u(x,t)$$

And

$$^cD_t^\alpha F(x,t) - F_{xx}(x,t) \geq 2\epsilon pu^{p-1}(x,t)u_x(x,t) + \epsilon p(p-1)(x-d_1)u^{p-2}(x,t)u_x^2(x,t).$$

Hence, the following inequality is obtained;

$$^cD_t^\alpha F(x,t) - F_{xx}(x,t) \geq 0 \text{ in } (d_1,d_2) \times [\tau,T).$$

*where $u_x(x,t) > 0$ and $p > 1$.*

Now, since $\epsilon$ is small positive enough, then the values at the boundary are obtained as follows;

$$F(d_1,t) = u_x(d_1,t) - \epsilon(d_1 - d_1) \cdot u^p(d_1,t) = u_x(d_1,t) > 0$$

And

$$F(d_2,t) = u_x(d_2,t) - \epsilon(d_2 - d_1) \cdot u^p(d_2,t) > 0$$

Similarly, the following inequality is obtained on initial datum;

$$F(x,0) = u_x(x,0) - \epsilon(x-d_1) \cdot u^p(x,0) > 0$$

With the help of the maximum principle, it is obtained that $F(x,t) \geq 0$ in $[d_1,d_2] \times [\tau,T]$. So, $u_x(x,t) \geq \epsilon(x-d_1)u^p(x,t)$ in $[d_1,d_2] \times [\tau,T]$. Integrating the last inequality with respect to x from $d_1$ to $d_2$, the following inequality is obtained;

$$u(d_1,t) \leq \left(\frac{2}{\epsilon(p-1)(d_2-d_1)^2}\right)^{\frac{1}{(p-1)}} < \infty$$

This means the solution $u(x,t)$ remains finite for all $x \in [0,1)$ and $t \in [0,T]$. Therefore, it does not blow up in $(0,1)$ and a blow-up can only occur only at $x = 1$.

## 5. Conclusion

This article is concerned with the blow-up solutions of the time fractional heat equation subject to a nonlinear Neumann boundary condition of power type. Under certain restricted conditions, it is proven that every positive solution achieves blowup in a finite time. Moreover, the blow-up can only be achieved at the boundary points $x = 1$.

In future studies, the theoretical results obtained in this study are planned to be supported by numerical approaches. Another open problem is to study the quenching version of the same problem.